\documentclass[11pt]{article}

\input{epsf.sty}

\usepackage{amsfonts,amssymb}

\newcommand{\Z}{\mathbb{Z}}

\usepackage{amsfonts,amssymb,amsmath}
\newtheorem{theorem}{Theorem}[section]

\newtheorem{proposition}[theorem]{Proposition}

\newtheorem{definition}[theorem]{Definition}

\newtheorem{example}[theorem]{Example}

\newtheorem{remark}[theorem]{Remark}

\begin{document}

\title{Tangle Embeddings and Quandle Cocycle Invariants}
\author{Kheira Ameur \and Mohamed Elhamdadi \and Tom Rose
\and Masahico Saito\footnote{Supported in part by NSF Grant DMS \#0301089,
 \#0603876.}
 \and Chad Smudde
}

\date{}

\maketitle

\vspace{-10mm}
\begin{center}
{Department of Mathematics and Statistics} \\
{University of South Florida}
\end{center}

\begin{abstract}
To study embeddings of tangles in knots, we use quandle cocycle invariants. Computations are carried out for the tables of knots and
tangles, to investigate which tangles may or may not embed in knots
in the tables.
\end{abstract}

\section{Introduction}

A {\it tangle} is a
pair $(B, A)$ where $A$ is a set of properly embedded arcs in a $3$-ball $B$.
A tangle will have  four end points throughout the article, unless otherwise specified.
A tangle $T$ is {\it embedded} in a link (or a knot)  $L$
if there is an embedded  ball $B$ in $3$-space such that $T$
is equivalent to the pair $(B, B \cap L)$. All maps are assumed to be smooth.
Tangles are represented by diagrams in a manner similar to knot diagrams.

Tangle embeddings have been studied by several authors recently.
In \cite{Kre99}, the determinant was used in relation to evaluations of the
Jones polynomial, that have been further investigated
in \cite{CL05*,Kre99,KSW00}. Topological interpretations were
considered in \cite{PSW04*,Rub00}.
Tangles were also used to study
 DNA
recombinations~\cite{ES}.

In this article,  we present a method of using quandle cocycle
invariants as obstructions to embedding tangles in knots, and
examine their effectiveness as obstructions by looking at the table
of tangles presented in \cite{KSS03} and the knot table
in~\cite{Chuck}. Quandles are self-distributive sets with additional
properties (see below for details), and have been used in
the
study of knots since 1980s. A cohomology theory of quandles have
been developed, and their cocycles have been used as state-sum
invariants of knots and knotted surfaces~\cite{CJKLS}.
Quandles are  used to investigate
 tangles also  in
\cite{DN,maciej}.

In this paper, we focus on
effectiveness of quandle cocycle invariants as obstructions.
 It
will be shown that the invariants often provide effective
obstructions when a given tangle has non-trivial colorings by
quandles.

The paper is organized as follows. After a review %ME changed "expositions" to "review"
 of
preliminary material in Section~\ref{prelimsec}, colorings of
tangles are defined in Section~\ref{tanglecolsec}, and the tangles
in the table (see Fig.~\ref{TTable}) that have non-trivial colorings
by Alexander quandles are listed. The main theorem is presented in
Section~\ref{quandleobssec}. For tangles listed in
Section~\ref{tanglecolsec}, it is examined which tangles may or may
not embed in knots in the knot table. In
Section~\ref{manytanglesec}, embeddings of multiple disjoint copies
of tangles are discussed. Part of the results are based on the work
in the Ph.D. dissertation by Kheira Ameur \cite{Kheira}.

\section{Preliminaries}\label{prelimsec}

\subsection{Tangles and their operations}

The  conventions described in this subsection
are commonly found in the literature (see for example, \cite{Ad94,Mura}).

The four end points of a given tangle diagram $T$
are located at four corners of a circle
in a plane at angles
$\pi/4$, $3\pi /4$, $5 \pi / 4$ and $7 \pi / 4$, when the circle is placed
with the origin as its center.
These end points
are labeled by NE, NW, SW, and SE, respectively,
representing North East, etc.

\begin{figure}[htb]
\begin{center}
\mbox{
\epsfxsize=1.2in
\epsfbox{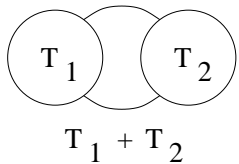}
}
\end{center}
\caption{ Addition of  tangles }
\label{addition}
\end{figure}

\begin{figure}[htb]
\begin{center}
\mbox{
\epsfxsize=1.7in
\epsfbox{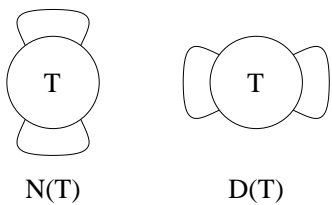}
}
\end{center}
\caption{ Closures (numerator $N(T)$ and
 denominator $D(T)$) of
tangles } \label{NtDt}
\end{figure}

The addition $T_1+T_2$ of two tangles $T_1$, $T_2$ is another tangle
defined from the original two as depicted in Fig.~\ref{addition}.
There are two ways of closing the end points of a tangle, called
closures, the {\it numerator} $N(T)$ and {\it denominator} $D(T)$ of
a tangle $T$, defined as depicted in Fig.~\ref{NtDt}.

\begin{figure}[htb]
\begin{center}
\mbox{
\epsfxsize=3in
\epsfbox{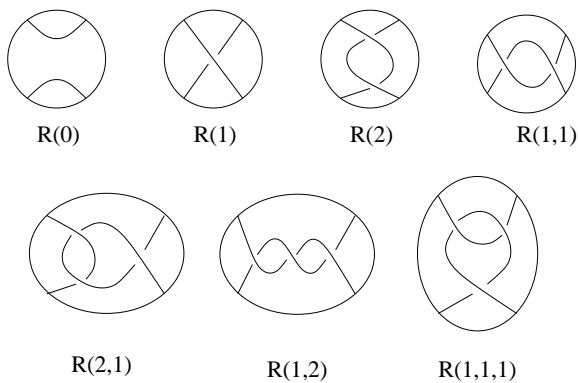}
}
\end{center}
\caption{ Some rational  tangles }
\label{rational}
\end{figure}

There is a family of ``trivial'' or ``rational''  tangles, some of
which are depicted in Fig.~\ref{rational}. These tangles
are obtained from the trivial tangle of two vertical straight arcs by
successively twisting end points vertically and horizontally. See
again \cite{Ad94} or \cite{Mura}   for more details.

\begin{figure}[htb]
\begin{center}
\mbox{
\epsfxsize=4.5in
\epsfbox{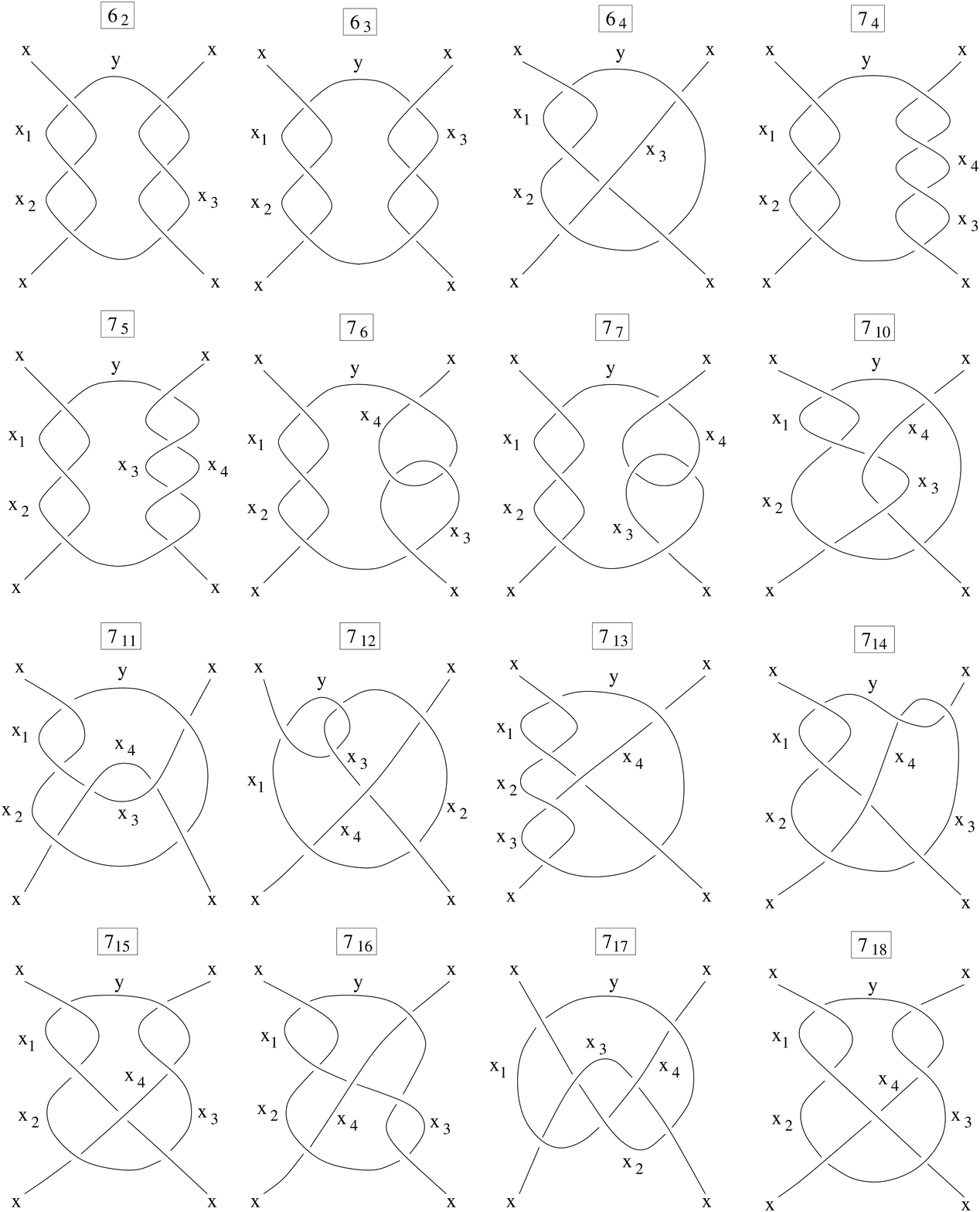}
}
\end{center}
\caption{ Tangles with non-trivial colorings by Alexander quandles}
\label{TTable}
\end{figure}

In \cite{KSS03},
 {\it prime} tangles (with crossing number at most seven)
  are classified,
and a
table of their diagrams is
given. The table consists of a single $5$ crossing tangle followed
by four $6$ crossing tangles, and $18$ tangles of $7$ crossings. Some
multiple component tangles were also classified. The tangles are
named in a scheme similar to
knots by integers
with subscripts.
Some of the tangles are presented in Fig.~\ref{TTable}.

\subsection{Quandles, colorings, and cocycle invariants}

A {\it quandle}, $X$, is a set with a binary operation
$(a, b) \mapsto a * b$
such that

(I) For any $a \in X$,
$a* a =a$,

(II) For any $a,b \in X$, there is a unique $c \in X$ such that
$a= c*b$,

(III)
For any $a,b,c \in X$, we have
$ (a*b)*c=(a*c)*(b*c). $

\noindent
A {\it rack} is a set with a binary operation that satisfies
(II) and (III).
Racks and quandles have been studied in, for example,
\cite{Br88,FR,Joyce,Matveev}.

The following are typical examples of quandles.
A group $G$
with
conjugation
as the quandle operation,
$a*b=b a b^{-1}$,
is a quandle.
Any ${\Z }[t, t^{-1}]$-module $M$
is a quandle with
$a*b=ta+(1-t)b$, $a,b \in M$,
that is
called an {\it  Alexander  quandle}.
Let $n$ be a positive integer, and
for elements  $i, j \in \Z_n$,  define $i\ast j \equiv 2j-i \pmod{n}$.
Then $\ast$ defines a quandle
structure  called the {\it dihedral quandle},
  $R_n$.

\begin{figure}
\begin{center}
\mbox{
\epsfxsize=2.5in
\epsfbox{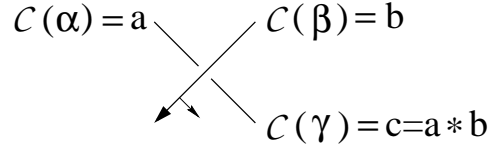}
}
\end{center}
\caption{ Quandle relation at a crossing  }
\label{qcolor}
\end{figure}

Let $X$ be a fixed quandle.
Let $K$ be a given oriented classical knot or link diagram,
and let ${\cal R}$ be the set of (over-)arcs.
The normals (normal vectors) are given in such a way that
the ordered pair
(tangent, normal) agrees with
the orientation of the plane, see Fig.~\ref{qcolor}.
A (quandle) {\it coloring} ${\cal C}$ is a map
${\cal C} : {\cal R} \rightarrow X$ such that at every crossing,
the relation depicted in Fig.~\ref{qcolor} holds.
Specifically, let $\beta$ be the over-arc at a crossing, and let
$\alpha $ and $\gamma$ be the under arcs, such that the normal of
the over-arc points from $\alpha$ to $\gamma$, then $\cal
C(\alpha)*\cal C(\beta)=\cal C(\gamma)$ holds.
The  (ordered) colors $( {\cal C}(\alpha), {\cal C}(\beta)) $
are called {\it source} colors.
Let ${\rm Col}_X(K)$ denote the set of colorings of a knot diagram $K$
by a quandle $X$.

Let $K$ be a knot diagram on the plane. Let $X$ be a finite quandle
and $A$ be an abelian group. Let  $\phi : X \times X \rightarrow A$
be a  quandle $2$-cocycle, which can be regarded as a function
satisfying the $2$-cocycle condition
$$\phi(x,y) - \phi(x,z) + \phi(x*y, z) - \phi(x*z, y*z)=0, \quad \forall x,y,z \in X$$
and $\phi(x,x)=0, \forall x \in X$.
Let ${\cal C}$ be a coloring of a given knot diagram $K$ by $X$.

The {\em Boltzmann weight} $B( {\cal C}, \tau)=B_{\phi}( {\cal C}, \tau)$
at a crossing $\tau$ of $K$
is then defined by  $B( {\cal C}, \tau)=\epsilon(\tau) \phi(x_{\tau}, y_{\tau})$,
where $(x_{\tau}$, $y_{\tau})$ is  the
source
 colors at $\tau$ and
$\epsilon(\tau)$ is the sign ($\pm 1$) of $\tau$.
Then the $2$-cocycle invariant $\Phi(K)=\Phi_{\phi}(K)$
in a multiset form is defined by
$$\Phi_{\phi}(K) =  \left\{ \left. \sum_{\tau} B( {\cal C}, \tau ) \ \right|
 \ {\cal C}\in {\rm Col_X(K)} \right\} .$$
 (A multiset
 is
 a collection of elements where a single element can be repeated
 multiple times, such as $\{ 0, 0, 1, 1, 1 \}$,
 which is also denoted by $\{ \sqcup_2 0 , \sqcup_3 1 \}$).

Let  $\theta : X \times X \times X \rightarrow A$
be a  quandle $3$-cocycle, which can be regarded as
a function satisfying
\begin{eqnarray*}
\lefteqn{\theta (x,z,w)- \theta(x,y,w)+\theta(x,y,z)
-\theta(x*y,z,w)} \\ & & +\theta(x*z,y*z,w)-\theta(x*w,y*w,z*w)=0,
 \quad \forall x,y,z, w \in X,
 \end{eqnarray*}
and $\theta(x,x,y)=0=\theta(x,y,y), \forall x, y \in X$.

\begin{figure} % [htb]
\begin{center}
\mbox{
\epsfxsize=2in
\epsfbox{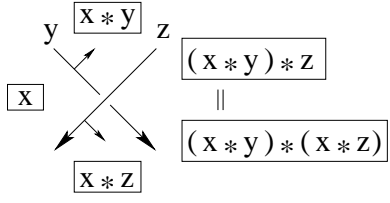}
}
\end{center}
\caption{ Region colors at a crossing  }
\label{regioncolors}
\end{figure}

Let ${\cal C}$ be a coloring of arcs and regions of a given diagram
$K$. Specifically, for a coloring $ {\mathcal C}$, there is a
coloring of regions that extend $ {\mathcal C}$ as depicted in
Fig.~\ref{regioncolors}. Suppose that two regions $R_1$ and $R_2$
are separated by  an arc colored by $y$, and the normal of the arc
points from $R_1$ to $R_2$. If $R_1$ is colored by $x$, then $R_2$
receives the color $x*y$. Let $(x_{\tau}, y_{\tau}, z_{\tau})$
(called the {\it  ordered triple of colors}  at a crossing  $\tau$ )
be the colors near a crossing $\tau$ such that $x$ is the color of
the region (called the source region) from which both orientation
normals of over- and under-arcs point, $y$ is the color of the
under-arc (called the source under-arc) from which the normal of the
over-arc points, and $z$ is the color of the over-arc. See
Fig.~\ref{regioncolors}. Let $(x_{\tau}, y_{\tau}, z_{\tau})$ be the
ordered triple of colors at a crossing  $\tau$. Then the weight in
this case is defined by
 $B( {\cal C}, \tau)=\epsilon(\tau) \phi(x_{\tau}, y_{\tau}, z_{\tau}).$
The $3$-cocycle invariant is defined in a similar way to the
$2$-cocycle invariant
 by the multiset $ \Phi_{\theta}(K) =\{
\sum_{\tau} B( {\cal C}, \tau ) \ | \ {\cal C}\in {\rm Colr_X(K)} \}
$, where  ${\rm Colr_X(K)}$ denotes the set of colorings with region
colors of $K$ by $X$.

If the quandle $X$ is finite, the invariant as a multiset can be
written by an expression similar to those for the state-sums;
  if a given
multiset of group elements is $\{ \sqcup_{m_1} g_1, \ldots,
\sqcup_{m_\ell} g_\ell \} $, then we use the polynomial notation
$m_1 u^{g_1} + \cdots + {m_\ell} u^{g_\ell }$ where
 $u$ is
 a formal symbol. For example, the multiset value of
the invariant for a trefoil with the Alexander quandle $X=\Z_2[t,
t^{-1}]/ (t^2+t+1)$ with the same coefficient $A=X$ and  a certain
$2$-cocycle is $\{ \sqcup_4 (1) ,  \sqcup_{12} (t+1) \}$, and
 is
denoted by $4 + 12 u^{(t+1)}$, where we use the convention $u^0=1$ and
exponential rules apply.

For computing the invariants, one needs an explicit formula for
cocycles. Polynomial expressions were used first
 in \cite{Mochi},
and investigated closely including  higher dimensional cocycles in
\cite{Kheira}.

\section{Boundary monochromatic colorings and the cocycle invariants of tangles}
\label{tanglecolsec}

We use quandle cocycle invariants as obstructions to embedding tangles
in knots.
We first define cocycle invariants for tangles.

\begin{definition}{\rm
\begin{sloppypar}
Let $T$ be a tangle and $X$ be a quandle. A (boundary-monochromatic)
{\it coloring} ${\cal C}: {\cal A} \rightarrow X$ is a map from the
set of arcs in a diagram of $T$ to $X$ satisfying the same quandle
coloring condition as for knot diagrams at each crossing, such that
the (four) boundary points of the tangle diagram receive the same
element of $X$.
\end{sloppypar}

For a coloring ${\cal C}$ of a tangle diagram $T$,
a region colorings are defined in a similar manner as  in the knot case.
In this case, we allow region colors to change (not necessarily
colored by the same element as the one assigned to the boundary points).
} \end{definition}

Denote by ${\rm Col}_x(T)$ and ${\rm Col}_X(T)$
the set of boundary-monochromatic colorings of $T$
with the boundary color $x \in X$
and the set of all   boundary-monochromatic colorings, respectively.
Let $\Phi(T,x)= \sum_{C\in
Col_x(T)}\prod_{\tau}B( C, \tau)$. Then the cocycle invariant for a tangle $T$
is defined by $
\Phi_{\phi}(T)=\sum_{x\in X}\Phi(T,x)$.
The invariants with region colors are defined in a similar manner,
by taking sum over all colorings of regions as well as  colorings of diagrams.

It is seen in a way similar to the knot case that the number of
colorings $|{\rm Col}_X(T)|$ does not depend on a choice of a
diagram of $T$. Any coloring ${\cal C}_1$ of a diagram of a tangle
$T$ is changed via a sequence of Reidemeister moves
(with boundary points fixed) %% added
 to a coloring
${\cal C}_2$ of $T$.  Given two diagrams $D_1$ and $D_2$ of a tangle
$T$, there is one-to-one correspondence between the set of
colorings of $D_1$ and the set of colorings of $D_2$ and the cocycle invariant
is well-defined.

\begin{table}[htb]

\begin{center}

\begin{tabular}{|l||l|}\hline
Quandle & Tangle colored\\ \hline $
\mathbb{Z}_p[t,t^{-1}]/(t^2-t+1)$ & $6_2$, $6_3$, $7_{17}$(NW In, SW
In). \\ \hline $\mathbb{Z}_2[t,t^{-1}]/(t^2+t+1)$ & $6_2$, $6_3$,
 $7_{17}$(NW In, SW In), \\
 & $7_4$(NW In, NE In),
$7_5$(NW In, NE In), \\
 & $7_6$(NW In, NE In), $7_7$(NW In, NE In), \\
 \hline $R_3$ & $6_2$, $6_3$, $7_{16}$,
$7_{17}$. \\ \hline $R_5$ & $7_{13}$, $7_{18}$.\\
\hline $R_7$ & $7_{15}$.\\
\hline

\end{tabular}

\end{center}

\caption{Tangles with non-trivial colorings}\label{QuTanCol}
\end{table}

Table~\ref{QuTanCol} summarizes the tangles in the tangle table
 that have non-trivial boundary monochromatic colorings
by some Alexander quandles.
These are found by hand calculations, occasionally assisted by {\it Maple}.
Specifically, variables $x_i$, $i=1,2,\ldots$, are assigned on arcs
of tangle diagrams as indicated in Fig.~\ref{TTable},
and the coloring conditions of the form $x_k=t x_i + (1-t) x_j$ are imposed
corresponding to crossings, and the system of linear equations
in $\Z[t, t^{-1}]$ are solved to find which Alexander quandles have
non-trivial colorings.
 These tangles
 with non-trivial colorings by Alexander quandles
 are depicted in Fig.~\ref{TTable}.
 This list  compares with the original list in \cite{KSS03} as follows.
 Their list starts with
 one $5$-crossing tangle $5_1$, which colors trivially by Alexander quandles
 and is not listed in Fig.~\ref{TTable}.
 There are four $(6_1 - 6_4)$ $6$-crossings tangles, three of which are
 in our list.
 Thirteen out of eighteen $7$-crossing tangles are in our list.

%\clearpage % adjust this

\section{Quandle cocycle invariants as obstructions to tangle embeddings}
\label{quandleobssec}

The quandle $2$- and $3$-cocycle invariants are defined for tangles in a
manner similar to the knot case
using the set of boundary monochromatic colorings, and denoted by $\Phi_{\phi}(T)$.
We use the multiset version of the invariant.

\begin{definition}{\rm
The inclusion of multisets is
denoted by $\subset_m$.
Specifically, if an element $x$ is repeated $n$ times in a multiset,
call $n$ the multiplicity of $x$, then $M \subset_m N$
for multisets $M$, $N$ means that
if $x \in M$, then $x \in N$ and the multiplicity of $x$ in $M$ is less than or equal
to the multiplicity of $x$ in $N$.
}\end{definition}

\begin{theorem}\label{tanglethm}
Let $T$ be a tangle and $X$ a quandle.
Suppose $T$ embeds in a link $L$.
Then we have the inclusion
$\Phi_{\phi}(T) \subset_m \Phi_{\phi}(L)$.
\end{theorem}
{\it Proof.\/}
Suppose a diagram of $T$ embeds in a diagram of $L$.
We continue to use $T$ and $L$ for these diagrams.
For a coloring ${\cal C}$ of $T$, let $x$ be the color of the boundary points.
Then there is a unique coloring ${\cal C}'$ of $L$ such that the restriction of ${\cal C}'$
on $T$ is ${\cal C}$ and all the arcs of $L$ outside of $T$ receive the color $x$.
 Then the contribution of $\sum_{\tau \in T} B( {\cal C}, \tau )$ to $\Phi_{\phi}(T)$
 is equal to the contribution $\sum_{\tau \in L} B( {\cal C}', \tau )$ to $\Phi_{\phi}(L)$,
 and the theorem follows.
 The same argument works for region colors and $3$-cocycle invariants.
 $\Box$

\clearpage % \bigskip

In Table~\ref{resultstable}, a summary is presented for the
tangles that color non-trivially by Alexander quandles.
 In the left
column of the table, the tangles that appear in Table~\ref{QuTanCol}
are listed. In the middle column, knots % in the table
that we found
to embed % in
a given tangle are listed. The third column lists the
knots for which we could not exclude the possibility of embedding of
the given tangle using
 cocycle invariants. The tangles are specified
by the notation $T(6_2)$, for example, for the tangle numbered $6_2$
to distinguish them from knots. We note that there are $84$ knots in
the table up to (including) $9$-crossing knots. For the tangle
$T(6_3)$, for example, all except $3$ out of $84$ are detected by
the cocycle invariants that they do not embed the tangle. It is
checked by hand that these remaining three do embed it.
% delete : %%%
%The cocycle invariant is least effective for the tangle  $T(7_6)$,
% in the table is for $7_6$,
%for which $23$ out of $84$ knots
%remain not excluded.
%This is partly because
%the Alexander quandles that can color this tangle
% non-trivially are limited.

\begin{table}
\begin{center}

\begin{tabular}{|l||l|l|}\hline
Tangle & Embeds in: & May embed in: \\ \hline \hline
  $6_2$ (NW In, SW Out) & $(8_5)^*=N(T(6_2)+R(-2)) $ & $8_{18}, 9_{29}, 9_{38}$.\\
  \hline
  $6_2$ (NW In, SW In) & $(3_1)=N(T(6_2)+R(-1))$
  & $3_1, 7_4, 7_7, 8_{18}, 9_{10}, 9_{29},$ \\
  & & $9_{35}, 9_{37}, 9_{38}, 9_{46}, 9_{48}.$\\ \hline
  $6_3$ (NW In, SW Out) & $(8_{10})=N(T(6_3)+R(2,1))^*$ &  \\
 &  $(8_{20})=N(T(6_3)+R(2))^* $&    $8_{10}, 8_{20}, 9_{24}$.\\
  & $(9_{24})=N(T(6_3)+R(2,2))^* $& \\  \hline
  $7_4$ (NW In, NE In) & $(4_1)=(N(T(7_4)+R(-1))$
  & $3_1,  4_1, 7_2, 7_3, 8_1, 8_4, 8_{11}, $ \\
  & & $ 8_{13}, 8_{18},  9_1, 9_6, 9_{12}, 9_{13}, $ \\
  & & $9_{14}, 9_{21}, 9_{23}, 9_{35}, 9_{37}, 9_{40}.$\\ \hline
  $7_5$ (NW In, NE In) & $(7_3)^*=N(T(7_5)+R(-1))$ &
    Same as  $7_4$(NW In, NE In).  \\ \hline
  $7_6$ (NW In, NE In) & &
  $ 8_5,
8_{10}, 8_{15}, 8_{18} - 8_{21},  9_{16},$\\
& & $  9_{22}, 9_{24}, 9_{25},
9_{28} - 9_{30}, 9_{36},  $ \\
& & $9_{38}, 9_{39}, 9_{41} - 9_{45}, 9_{49}$.\\ \hline
  $7_7$ (NW In, NE In) & & Same as  $7_6$(NW In, NE In).\\ \hline
  $7_{13}$ (NW In, NE Out) & $(7_4)=N( T(7_{13}) )$ & \\
  &  $(8_{16}) = N( T(7_{13}) + R(1) )$  &    \\

& $ (9_{39}) = N( T(7_{13}) + R(1,1) ) $ &
    $4_1, 7_4, 9_{24},  9_{37}, 9_{39}, 9_{40}, 9_{49}.$ \\
 & $ (9_{49}) = N( T(7_{13}) $ & \\
  & \mbox{\hspace{20mm}} $+ R(-1,-1) )$  & \\  \hline
  $7_{15}$ (NW In, SW In) & $ (5_2) =  N( T(7_{15}) + R(-1) )$&
  $ 5_2, 8_{16}, 9_{41}, 9_{42}$.  \\  \hline
  $7_{15}$ (NW In, SW Out) & $(7_7)=D( T(7_{15}) )$ &
     $7_1, 7_7, 8_5, 9_4, 9_{12}, 9_{41}$.\\
  &  $(9_{41})=N( T(7_{15}) + R(2))$& \\  \hline
  $7_{16}$ (NW In, NE In) &  $(7_7)^* = D( T(7_{16}) ) $
 &
 $8_5,8_{15}, 8_{18}, 8_{19}, 8_{21}, 9_2, 9_4,  $
 \\
& &   $ 9_{11},  9_{15}, 9_{16},9_{28}, 9_{34}, 9_{37},$ \\
& & $ 9_{40}, 9_{46}, 9_{47}$.
\\ \hline
 $7_{16}$ (NW In, NE Out) & $(7_4) = N( T(7_{16}) )$ &
 Same as $6_2$(NW In, SW In). \\ \hline
 $7_{17}$ (NW In, SW In) & $(8_{18} )=N( T(7_{17}) + R(1) ) $
 & $8_{18}, 9_{40}$. \\ \hline
  $7_{17}$ (NW In, SW Out) &
 & Same as $7_{16}$(NW In, NE In). \\ \hline
  $7_{18}$ (NW In, SW In) & $ (8_{21})=N( T(7_{18}) + R(1) ) $&
  $5_1, 8_{18}, 8_{21}, 9_2, 9_{12}, 9_{23}, 9_{31}, $ \\
   & & $9_{40}, 9_{49}$. \\ \hline
   $7_{18}$ (NW In, SW Out) & $ (5_1)= D( T(7_{18}) )$
 & Same as $7_{18}$(NW In, SW In). \\
\hline
\end{tabular}
\end{center}

\caption{Summary of the
results}\label{resultstable}
\end{table}

To demonstrate how we obtain
 these results, we state and prove
  the following.

\begin{proposition}\label{prop62}
The tangle $T(6_2)$ with the orientation of the NW arc inward and
the SW arc outward does not embed in the  knots in the table up to
$9$ crossings  except, possibly,  for  $8_{18,}$, $ 9_{29}$,
$9_{38}$.
\end{proposition}
{\it Proof.\/}
 With this orientation,  the tangle is of the form of two copies of  the mirror
 of the trefoil, and is colored non-trivially by   the quandle
 $Z_p[t]/(t^2-t+1)$.

\begin{sloppypar}
We exhibit a method to determine the invariant from the
table in \cite{Chad}. For $p=2$, the table of quandle cocycle
invariants in \cite{Chad} gives $16 + 48u^t$ as the invariant for
trefoil
with the $3$-cocycle $\phi(x,y,z)=(x-y)(y-z)^2$.
This implies that any non-trivial coloring contributes $t$
 to the invariant. Its mirror has the same property.
%% Delete :
% (Note that this case p=2 gives the same values of the invariant for mirror images.)
 With two copies, any non-trivial coloring of the tangle  contributes $2t=0$ when $p=2$.  Hence the invariant value of the tangle is $64$.
 From the table this does not embed in knots up to $9$ crossings except for the following possibilities :
 $8_5$, $8_{10}$, $8_{15}$, $8_{18}$, $8_{19}$, $8_{20}, 8_{21}, 9_{16}, 9_{22}, 9_{24}, 9_{25}, 9_{28}, 9_{29}, 9_{30},  9_{36}, 9_{38}, 9_{39}, 9_{40}, 9_{41},
 9_{42}$,
$ 9_{43}$,
 $9_{44}$,
  $9_{45}, 9_{49}$.
  \end{sloppypar}

 For $p=3$,  the invariant table gives $243 + 486u^{(2t+2)}$
  as the invariant for trefoil. This implies that $486$ non-trivial colorings contributes
  $2t+2$ to the invariant. Its mirror contributes $t+1$. With two copies, $486$
   non-trivial colorings of the tangle  contributes $2t+2$.
   Hence the invariant value of the tangle is $243 + 486u^{(2t+2)}$.
   {}From the table this does not embed in knots up to $9$ crossings except for:
   $3_1, 8_{18}, 9_2, 9_4, 9_{29}, 9_{34}, 9_{38}$.

 For $p=5$, the table gives  $$625+3750u^{(t+3)}+3750u^{(4t+2)}+3750u^{(3t+4)}+3750u^{(2t+1)} $$
 as the invariant for the trefoil. As in the previous cases, the tangle has the invariant value  $$625+3750u^{(3t+4)}+3750u^{(2t+1)}+3750u^{(4t+2)}+3750u^{(t+3)}$$
  (for example, for the contribution $t+3$ of trefoil, the mirror contributes $4t+2$,
   its double contributes $3t+4$).
   % D: % This is the same as the  trefoil.
   {}From the
   table this does not embed in knots up to 9 crossings except for:
   $3_1, 8_3, 8_5, 8_{11}, 8_{15}, 8_{18}$, $8_{19}, 8_{21}, 9_1, 9_5, 9_6, 9_{16}, 9_{19},$ $9_{23}, 9_{28}, 9_{29}, 9_{38}, 9_{40}$.
 %%   (by the symmetry of the  invariant values).

 For $p=7$, the trefoil has $117649$
  as the invariant value, and so does the tangle. From the
  table this does not embed in knots up to 9 crossings except for:
  $3_1, 8_5, 8_{10}, 8_{11}, 8_{15}, 8_{18}$, $8_{19}, 8_{20}, 8_{21}, 9_1,  9_6, 9_{16},  9_{23}, 9_{28}, 9_{29}, 9_{38}, 9_{40}$.

 From all these information combined,
  this tangle does not embed in knots up to $9$ crossings except
 for the only possibilities of
 $ 8_{18},  9_{29},  9_{38}$.
   $\Box$

\bigskip

We have not been able to determine whether the tangle $T(6_2)$
actually embeds in these three knots that the invariant failed to
exclude. In the next example, however, we were able to determine
completely the embedding problem up to $9$ crossings.

\begin{proposition}\label{prop63}
The knots in the table up to $9$ crossings in which the
 tangle $T(6_3)$ embeds are
 exactly  $8_{10}$, $8_{20}$, $9_{24}$.
Here, the orientation of the tangle is such that
  the end point NW is oriented inward
 and the SW end point is oriented  outward.
\end{proposition}
{\it Proof.\/}
 The tangle $T(6_3)$ is written as the addition
 $R(3)+R(-3)$. Hence it is
  colored   non-trivially by $Z_p[t]/(t^2-t+1)$
 for any $p \in \Z$ (we use only primes), as well as the dihedral quandle $R_3$.
  For the quandle $Z_p[t]/(t^2-t+1)$ we used the $3$-cocycle
 $f(x,y,z)=(x-y)(y-z)^p$.
The colors of the source region for these two copies of the trefoil diagrams
 ($R(3)$ and $R(-3)$)  coincide.  The signs of the crossings are opposite. Hence the invariant is trivial, $(p^2)^3$ copies of $0$, for $Z_p[t]/(t^2-t+1)$.
   For $p=5$, in particular,    from the calculations
    in \cite{Chad},
    Theorem  \ref{tanglethm}     implies that this tangle
    may  embed, among knots in the table up to $9$ crossings,
     only in: $8_{10}$, $8_{12}$, $8_{18}$, $8_{20}$, $9_{24}$.
The invariant with $R_3$ further excludes $8_{12}$ and $8_{18}$.
Therefore the tangle may embed only in $8_{10}$, $8_{20}$, and
$9_{24}$.

On the other hand, it is seen that
\begin{eqnarray*}
 (8_{10}) &=& N( T(6_3) + R(2,1) )^*, \\
 (8_{20}) &=& N( T(6_3) + R(2) )^*, \\
 (9_{24}) &=& N( T(6_3) + R(2,2) )^* ,
\end{eqnarray*}
where $K^*$ denotes the mirror image of a knot $K$, and $R$ denotes
the rational tangles. Note that this tangle $T(6_3)$ is equivalent
to its mirror. Therefore we have shown that the tangle $T(6_3)$ does
indeed embed in these three knots. $\Box$

\begin{remark}{\rm
In general the orientation needs to be specified to define the
quandle cocycle invariants. (In our case only the dihedral quandles
can be used for the invariant without specifying the
orientations~\cite{ShinGraph}.) Furthermore, the mirror images of a
given knot in the table may be different. Thus all of our results
are stated for oriented tangles and oriented knots, and do not
include their mirror images. Our convention for specifying
orientations of tangles are already explained. For knots in the
table, we used Livingston's table~\cite{Chuck}, which includes
particular choices of mirrors if a knot is not amphicheiral. For the
orientations, we used the braid form in \cite{Chuck}, for our
calculations, so that the orientations are specified by downward
orientations of the braids.

%% Delete (may confuse that the mirror has a formula) %%%
%In Table~\ref{resultstable}, a few specific embeddings into
%mirrors of knots are mentioned. This information is not relevant
%because of the remark in the preceding paragraph, but these facts
%were used in computing the quandle cocycle invariants
%using the invariant values of the mirror of the knot in question.
%This method is indicated in the proof of Proposition~\ref{prop62}.

} \end{remark}

\section{Embedding disjoint tangles}\label{manytanglesec}

In this section we discuss
embeddings of
disjoint union of
tangles in knots. We  prove a theorem that will be used as
obstruction to embedding disjoint union of tangles and  give some
examples.

Let $C=\sum_{i=1}^k m_i u^{c_i}$, $ D=\sum_{j=1}^\ell  n_j u^{d_j}$
be polynomial expressions of multisets values of the invariants,
where
$m_i, n_j \in \Z_+$, $c_i, d_j \in A$,
where $A$ is  the coefficient abelian group.
Then we define $C \times D=\sum_{i, j} m_i n_j u^{c_i + d_j}$.
Let $|X|$ denote the number of elements of a quandle $X$.

% Move down %%%%
% Let $\phi$ be a $2$-cocycle of a quandle $X$ and we
%  consider $2$-cocycle invariants.
The quandle cocycle invariants are defined
for disjoint union of tangles $T_1 \sqcup \cdots \sqcup T_k$ in a
manner similar to tangles requiring that all the boundary points of
$T_1, \ldots, T_k$ receive the same color.
Let $\phi$ be a $2$-cocycle of a quandle $X$ and we define % moved down here
% Thus  we define
$$\displaystyle \Phi_{\phi}(T_1 \sqcup \cdots \sqcup
T_k)=\sum_{x_j\in X}\prod_{i=1}^{k}\Phi_\phi(T_i,x_j) .$$

\begin{proposition}\label{DisTangthm}
Let $\phi$ be a $2$-cocycle.
Let $T_1, \ldots, T_k$ be a disjoint union of tangles
such that for all $i=1, \ldots, k$, the condition
 $\Phi_\phi(T_i,x)=\Phi_\phi(T_i,y)$ holds for all $x,y \in X$.
Then we have
$$\Phi_{\phi}(T_1 \sqcup \cdots \sqcup T_k)=\frac{1}{|X|^{k-1}}\Phi_{\phi}( T_1 )\times\cdots
\times\Phi_{\phi}( T_k ).$$
 Furthermore  if  a disjoint union of $T_1,\ldots, T_k$
embed in a link $L$, then
$$ \Phi_{\phi}(T_1 \sqcup \cdots \sqcup T_k)
\subset_{m} \Phi_{\phi}(L) . $$
\end{proposition}
{\it Proof. \/} We compute
  $$\displaystyle \Phi_{\phi}(T_1 \sqcup \cdots \sqcup
T_k)=\sum_{x_j\in
X}\prod_{i=1}^{k}\Phi(T_i,x_j)
=|X| \prod_{i=1}^{k}\Phi(T_i,x) $$ for any fixed $x \in X$ since
$\displaystyle \Phi(T_i,x)=\Phi(T_i,y)$ for all $x,y \in X$. The
condition also implies that $\displaystyle
\Phi(T_i,x)=\frac{1}{|X|}\Phi_\phi(T_i)$ for all $i=1, \ldots,k$.
Hence
 $$\Phi_{\phi}(T_1 \sqcup \cdots \sqcup
T_k)=\frac{1}{|X|^{k-1}}\Phi_{\phi}( T_1 )\times\cdots
\times\Phi_{\phi}( T_k ).$$ Thus by the same argument as
the proof of Theorem~\ref{tanglethm}, if $T_1
\sqcup \cdots \sqcup T_k$
 embeds in a link $L$, we have
$$\frac{1}{|X|^{k-1}}\Phi_{\phi}( T_1 )\times \Phi_{\phi}( T_2
)\times \cdots \times \Phi_{\phi}( T_k )\
 \subset_{m} \Phi_{\phi}(L).  \quad  \Box$$

 \begin{example}\label{disjtang2ex}
 {\rm

 For the following examples, we consider the quandle
 $X=\mathbb{Z}_2[t,t^{-1}]/(t^2+t+1)$, and the 2-cocycle $f(x,y)=(x-y)^2y$.
The invariant values for this quandle are available in \cite{Chad}
(here we used knots up to $9$-crossings).
It is seen that the following tangles satisfy the condition required in
Proposition~\ref{DisTangthm}
 by direct calculations.
 Alternatively,  either of the triviality of the invariant, or the property
 that only the trivial colorings make trivial contributions to the invariant,
 implies the condition required.

\bigskip

\noindent {\bf (a)} We compute $\Phi_f(T(6_2)\sqcup
T(6_2))=\frac{1}{4}(16\times 16)=64$ (Proposition~\ref{prop62}).
Using~\cite{Chad}
we compare this invariant to the cocycle
invariant of knots in the knot table, and conclude
  that
$T(6_2)\sqcup T(6_2)$ does not embed in any knot in the knot table
up to 9 crossings.

 The invariant value of $T(6_3)$ is $\Phi_f(T(6_3)) =16$
by an argument similar to those used for $3$-cocycles
in the proof of Proposition~\ref{prop63}.
By theorem~\ref{DisTangthm} $\Phi_f(T(6_3)\sqcup
 T(6_3))=\frac{1}{4}(16\times 16)=64$. Hence
   $T(6_3)\sqcup T(6_3)$ does not embed
in any knot in the knot table up to 9 crossings.

The disjoint union $T(6_2)\sqcup T(6_3)$ also has the same invariant value $64$,
hence the same conclusion holds.

\bigskip

\noindent {\bf (b)} The invariant of the tangle $T(7_5)$ with
orientation (NW In, NE In) is
 $\Phi_f(T(7_5))=4+12u^{(t+1)}$. This can be seen from the fact that
 $T(7_5)$ embeds in the knot $(7_3)^*$ and the number of colorings
 by this quandle is
   the same for $T(7_5)$ and $(7_3)^*$,
 so that by an argument similar to the proof of Proposition~\ref{prop62},
 the tangle has the same invariant value as $(7_3)^*$ (see \cite{Chad}).
    Hence by Theorem~\ref{DisTangthm}
 we obtain
 $$\Phi_f(T(7_5)\sqcup T(7_5))=\frac{1}{4}(4+12u^{(t+1)})^2=40+24u^{(t+1)}.$$
 %% ( 1 + 3 u^{(t+1)} )^2 = 1 + 6 u^{(t+1)} + 9 u^{ 2(t+1) } , 2(t+1)=0 mod 2. %%
 Using ~\cite{Chad},
we compare this invariant to the cocycle invariant of knots in the
table, and we conclude  that  $T(7_5)\sqcup T(7_5)$ does not
embed in any knot in the table up to 9 crossings.

\bigskip

 \noindent {\bf (c)} Again by
Theorem~\ref{DisTangthm},
$$\Phi_f(T(6_2)\sqcup
 T(7_5))=\frac{1}{4}16(4+12u^{(t+1)})=16+48u^{(t+1)}.$$
  We find that  $T(6_2)\sqcup T(7_5)$ does not embed
in any knot in the knot table up to 9 crossings with possible exceptions of $8_{18}$
and $9_{40}$.

Since the invariant value for $T(6_3)\sqcup T(7_5)$ is the same
as $T(6_2)\sqcup T(7_5)$, we
obtain the same conclusion.

} \end{example}

Let $\psi$ be a $3$-cocycle of a quandle $X$ with coefficient group $A$.
Denote by $\Phi_\psi(T, x,  s)$ the $3$-cocycle invariant with the boundary color
$x \in X$ and the color of the leftmost region $s \in X$.
Then the $3$-cocycle invariant for disjoint union of tangles
$\sqcup_{i=1} ^k T_i$ is defined if $T_i$ satisfy the condition
$\Phi_\psi(T_i, x,  s)=\Phi_\psi(T_i, x',  s')$ for all $x, x', s, s' \in X$ for all $i=1, \ldots, k$,
and is defined in this case by
$$\displaystyle \Phi_{\psi}(T_1 \sqcup \cdots \sqcup
T_k)=\sum_{x_j\in X, s \in X}\prod_{i=1}^{k}\Phi_\phi(T_i,x_j, s)
= |X| \sum_{x_j\in X}\prod_{i=1}^{k}\Phi_\phi(T_i,x_j, s)
$$
for a fixed $s \in X$ (and the invariant does not depend on  this choice of $s$).
Note that the invariant does not depend on a fixed region color $s$
because of the above assumption.
Then the same argument as the proof of the preceding theorem can be
applied to show the following.

\begin{proposition}\label{ShadowTangthm}
Let $\psi$ be a $3$-cocycle.
Let $T_1, \ldots, T_k$ be a disjoint union of tangles
such that for all $i=1, \ldots, k$, the condition
 $\Phi_\psi(T_i,x, s)=\Phi_\psi(T_i,x', s')$ holds for all $x  , x', s, s' \in X$.
Then we have
$$\Phi_{\psi}(T_1 \sqcup \cdots \sqcup T_k)=\frac{1}{|X|^{2(k-1)}}\Phi_{\psi}( T_1 )\times\cdots
\times\Phi_{\psi}( T_k ).$$
 Furthermore  if  a disjoint union of $T_1,\ldots, T_k$
embed in a link $L$, then
$$ \Phi_{\psi}(T_1 \sqcup \cdots \sqcup T_k)
\subset_{m} \Phi_{\psi}(L) . $$
\end{proposition}

 \begin{example}  {\rm
 We used dihedral quandles $R_p$ with Mochizuki's cocycle~\cite{Mochi}
 $$\psi(x,y,z)=(1/p) (x-y) [ \ (2 z^p -y^p) - (2z-y)^p \ ] \pmod{p}. $$
 The invariant values are available in \cite{Chad} up to $12$-crossing knots
 for $p=3, 5$.
By  arguments similar to those in Example~\ref{disjtang2ex},
it is seen that the following tangles satisfy the condition required in
Proposition~\ref{ShadowTangthm}.

\bigskip

\noindent {\bf (a)} For the dihedral quandle $R_3$ and the
Mochizuki's $3$-cocycle $\psi$, the tangle $T(6_2)$ satisfies the
condition in Theorem~\ref{ShadowTangthm}. This is because $T(6_2)$
is the sum of two copies of part of trefoil diagrams, and the
trefoil has the property that any non-trivial coloring gives the
same non-trivial contribution to the cocycle invariant. Since
$\Phi_{\psi}(T(6_2))=9(1+2u)$,  by Theorem~\ref{DisTangthm} we
obtain $\Phi_{\psi} (T(6_2)\sqcup T(6_2)
)=\frac{1}{3^2}81(1+4u+4u^2)=9+36u+36u^2$. Using ~\cite{Chad} we
compare this invariant to the cocycle invariant of knots in the knot
table, and we find that  $T(6_2)\sqcup T(6_2)$ does not embed in any
knot in the knot table up to 11 crossings (there are $801$) except,
possibly, $8_{18}$ and $11a_{314}$. {}From the invariant value, the
number of colorings of $T(6_2)\sqcup T(6_2)$ is $81$, and among 801
knots in the table up to $11$ crossings, there are $40$ with at
least $81$ colorings. Hence the number of colorings alone can
exclude all but $40$ knots, but  the cocycle invariant is able to
exclude all but $2$.

\bigskip

\noindent {\bf (b)} $\Phi_{\psi}(T(6_3) \sqcup T(6_3))=\frac{1}{3^2}
3^3 3^3=81$ with $R_3$, and $T(6_3)\sqcup T(6_3)$ does not embed in
any knot in the knot table up to 11 crossings,  except possibly
$10_{99}$, hence the cocycle invariant excludes all $801$ knots but
one.

\bigskip

\noindent
{\bf (c)}
For the quandle $R_5$  and  tangles $T(7_{13})$ and $T(7_{18})$ both have
invariant $25(1+2u+2u^3)$, hence $T(7_{13}) \sqcup T(7_{18})$ has invariant
$25(5+u+4u^2+2u^3+4u^4)$. Thus it does not embed
 in any knot in the knot table up to 11 crossings,  except possibly
 $10_{103}$, $10_{155}$, $11a_{317}$, $11n_{148}$.

}\end{example}

\noindent
{\bf Acknowledgments.}
We would like  to thank J.S. Carter and S. Satoh for valuable conversations.

%%%%%%%%%%%%%%%%%%%%%%%%%%%%%%%%

\end{document}